\documentclass[preprint,12pt]{elsarticle}

\usepackage{amssymb}
\usepackage{amsmath}

\usepackage{hyperref}
\usepackage{cleveref}
\usepackage{booktabs}
\usepackage{tabularx}
\usepackage[numbers]{natbib}
\usepackage{url}

\journal{Engineering Applications of Artificial Intelligence}

\begin{document}

\begin{frontmatter}


\author{Donya Dabiri\fnref{label1}}
\ead{dabiri.d@northeastern.edu}

\author{Joshua DaRosa\fnref{label1}}
\ead{darosa.jo@northeastern.edu}

\author{Milad Saadat\fnref{label1}}
\ead{saadat.m@northeastern.edu}

\author{Deepak Mangal\fnref{label1}}
\ead{d.mangal@northeastern.edu}

\author{Safa Jamali\corref{cor1}\fnref{label1}}
\ead{s.jamali@northeastern.edu}
\cortext[cor1]{Corresponding author}

\title{A detailed and comprehensive account of fractional Physics-Informed Neural Networks: From implementation to efficiency}

\affiliation[label1]{organization={Department of Mechanical and Industrial Engineering, Northeastern University},
            addressline={360 Huntington Avenue},
            city={Boston},
            postcode={02115},
            state={Massachusetts},
            country={USA}}



\begin{abstract}
Fractional differential equations are powerful mathematical descriptors for intricate physical phenomena in a compact form. However, compared to integer ordinary or partial differential equations, solving fractional differential equations can be challenging considering the intricate details involved in their numerical solutions. Robust data-driven solutions hence can be of great interest for solving fractional differential equations. In the recent years, fractional physics-informed neural network has appeared as a platform for solving fractional differential equations and till now, efforts have been made to improve its performance. In this work, we present a fully detailed interrogation of fractional physics-informed neural networks with different foundations to solve different categories of fractional differential equations: fractional ordinary differntial equation, as well as two and three dimensional fractional partial differential equations. These equations are solved employing two numerical methods based on the Caputo formalism. We show that these platforms are generally able to accurately solve the equations with minor discrepancies at initial times. Nonetheless, since in Caputo formalism, the value of a fractional derivative at each point requires the function's value in all of its previous history, it is computationally burdensome. Here, we discuss strategies to improve accuracy of fractional physics-informed neural networks solutions without imposing heavy computational costs.
\end{abstract}




\end{frontmatter}



\section{Introduction}
\label{sec1}
The indispensable and widespread role of differential equations in depicting physical phenomena is indisputable \cite{braun1983differential}. In science and engineering generally one makes an effort to describe the dynamics of a system as flux of one variable while other variables change. Hence, differential equations are everpresent mathematical descriptors of physical systems across different disciplines such as thermofluid sciences \cite{holman1986heat,white1990fluid,cengel2008introduction}, chemistry \cite{gavalas2013nonlinear}, biology \cite{jones2009differential}, economy \cite{gandolfo1971economic, zhang2005differential}, medicine \cite{hoppensteadt2012modeling}, and more.
Naturally, more complex physical phenomena require additional parameters and different derivatives within their descriptions. Fractional differential equations (FDEs) present a class of differential equations that are extremely powerful in making complicated descriptions compact and concise. From biology \cite{MAGIN20101586}, to electrochemistry \cite{OLDHAM20109}, economics \cite{tarasov2019history}, rheology \cite{SCOTTBLAIR194721}, control of dynamic systems \cite{caponetto2010fractional}, and even modeling COVID-19 outbreak \cite{shaikh2020mathematical}, FDEs have shown great promise in describing the system of interest in an efficient fashion. Nonetheless, numerically solving FDEs is far more challenging compared to integer order ordinary or partial differential equations (ODEs and PDEs), resulting in popularization of FDEs and their further widespread adoption for different applications.

Throughout the years, various methods have been proposed for solving FDEs and can be categorized into three groups of analytical, numerical and semi-analytical schemes. Laplace transform, Fourier transform, and Green's function are among the analytical methods that can be implemented to solve FDEs analytically \cite{podlubny1998fractional}. However, these methods are limited to relatively simpler dynamics. Semi-analytical methods, such as the homotopy analysis \cite{odibat2010reliable, zurigat2010homotopy} and Adomian decomposition \cite{momani2006decomposition}, offer enhanced accuracy by integrating analytical simplifications with numerical solutions, and are particularly effective in complex domains with the memory effects inherent in fractional derivatives. Despite their benefits, these methods may involve intricate setups and significant computational overhead compared to more straightforward numerical approaches. This complexity often prompts a transition to more direct and computationally manageable approaches. Numerical methods involve spectral methods \cite{zayernouri2014fractional}, finite element methods (FEM) \cite{zhao2017adaptive}, and finite difference methods (FDM) \cite{diethelm2002analysis}. In FDMs, the fractional derivative term is descritized based on the mathematical definition of slope, which is the first order derivative. Successively, higher-order derivatives are computed using chain rule in differentiation. 

Inspired by the common integer order derivative definition, Grünwald \cite{grunwald1867begrenzte} and Letnikov \cite{letnikov1868theory} proposed a straightforward definition for fractional derivative known as Grünwald–Letnikov (GL) method, which requires knowledge of a function over the whole interval ($-\infty$,\textit{t}] in order to compute fractional derivative of the function $f(t)$ at point \textit{t}. GL method is convergent for a limited range of functions, e.g. the functions that are either bounded in ($-\infty$,\textit{t}] or do not increase rapidly as $t \rightarrow -\infty$ \cite{garrappa2019evaluation}. To address the disadvantages corresponding to GL method adjustments are required. Firstly, a starting point $t_0$ is chosen for functions with unknown or undefined behavior in $(-\infty, t_0)$. With this approach, $f(t)$ should be approximated in $(-\infty, t_0)$ since the knowledge of the function is required in $(-\infty, t_0)$. Hence, $f(t)$ is approximated with the Taylor polynomial of $f$ centered at $t_0$ \cite{garrappa2019evaluation}. The order of polynomial would be the $m-1$ where $m$ is smallest integer higher than the fractional derivative order $\alpha$. This approach is the foundation of the Caputo method in computing fractional derivatives. With this approach, the unknown behavior of the function in $(-\infty, t_0)$ is resolved without causing any discontinuity at $t_0$ due to the proper order of polynomial approximation. Despite the advantages of Caputo over GL, its implementation poses some challenges. Considering the history-dependent feature and rather intricate numerical procedure involved in implementation of Caputo method for computing fractional derivatives, methods that facilitate these procedures and can provide robust and accurate solutions to FDEs can be transformative in their adoption to new physical systems.

In recent years, machine learning platforms have shown great promise in solving differential equations using various approaches, including learning the discretizations for PDEs \cite{bar2019learning}, implementing spectral methods in neural networks \cite{meuris2023machine}, or breaking the derivatives of a hidden state in an ODE to multiple steps and use hidden layers to compute derivatives in each step \cite{chen2018neural}. Physics-informed neural networks (PINNs) have emerged as a platform for solving forward and backward differential equations by training a neural network with respect to the governing equations, boundary and initial conditions \cite{PINNs_Raissi,raissi2020hidden,karniadakis2021physics}. Encoding the governing equations in training has made PINNs applicable to various fields, including fluid mechanics \cite{cai2021}, solid mechanics \cite{Haghighat2021}, rheology \cite{Mahmoudabadbozchelou2021,mahmoudabadbozchelou2021data, Saadat2022}, seismology \cite{rasht2022physics} and bioengineering \cite{herrero2022ep}, to mention a few. Seminal work of Pang and coworkers introduced application of PINNs to fractional derivatives, referred to as fPINNs, for solving forward and inverse advection-diffusion equations with the GL method \cite{pang2019fpinns}. Despite presenting a great potential in solving fPDEs, convergence is not always guaranteed in fPINNs. In a recent study, fPINNs were modified to a set of rheologically-relevant equations of interest (viscoelastic constitutive models), and Caputo method was employed to solve fPDEs in an inverse problem \cite{Dabiri2023}. Furthermore, fPINNs have also been employed alongside some adjustments to account for uncertainty quantification, known as stochastic fractional partial differential equations (SFPDEs) \cite{ma2023bi}. 

In this study, our objective is to solve forward FDEs by employing Caputo based methods, as the Caputo formalism has demonstrated to be more reliant for problems with initial conditions. We fully investigate the trade-off between accuracy and computational cost, and propose feasible methods to enhance accuracy while imposing the minimum burden on computational cost. The subsequent sections of the paper are organized as follows: In \cref{sec: method}, we first present a series of fDEs [of interest] to be solved alongside their initial/boundary conditions. Next, their corresponding analytical solutions are provided as the ground truth for benchmarking. We then briefly discuss the numerical methods for solving fractional derivatives. Finally, data-driven implementation of fDEs in neural networks is elaborated in details. In \cref{sec: result}, applicability of our fPINN approach to solving each set of fDEs is presented and discussed. Finally, in \cref{sec:conclusion}, a conclusion is given based on the preceding discussions.

\section{Problem setup and Methodology}\label{sec: method}
\subsection{Fractional equations}\label{frac_diff}
Since the ultimate goal of this work is to provide a generic platform for reliable and robust solution of different fractional differential equations through a data-driven approach, it is important to first define the set of equations to be evaluated and the benchmarking methods that are used as the ground truth for each test. The fDEs aimed for solving with the neural network are introduced in the order of complexity, starting with a fractional ordinary differential equation (fODE) and proceeding with fractional partial differential equations (fPDEs).\\

An example of a \textbf{Fractional Ordinary Differential Equation (fODE)} can be generally defined and written as:
\begin{equation}
    \frac{\partial^\alpha u(t)}{\partial t^\alpha} = -u(t) + t^2 + \frac{8}{3\sqrt{\pi}}t^{\frac{3}{2}},  \quad 0\leq t\leq 1, \alpha = 0.5
    \label{eq: fODE}
\end{equation}
Where $\alpha$ represents the fractional derivative order. We begin by simplifying the initial condition of $u(0) = 0$. The analytical solution to equation \ref{eq: fODE} then will be derived as:
\begin{equation}
    u(t) = t^2
    \label{eq: fODE_analytical}
\end{equation}

For \textbf{Fractional Partial Differential Equation (fPDE)}, we present the equations with respect to their level of complexity. For a \textit{\textbf{two-dimensional case}}, a fractional diffusion equation with a temporal term of fractional derivative order in two dimensions of $(x,t)$ is defined as follows:
\begin{equation}
    \frac{\partial^\alpha u(x,t)}{\partial t^\alpha} + u(x,t) = \frac{\partial^2 u (x,t)}{\partial x^2} + f(x,t) ,  \quad 0\leq x\leq 2, 0\leq t, 0<\alpha<1
    \label{eq: 2DfPDE}
\end{equation}
With the initial and boundary conditions defined as:
\begin{equation}
    \begin{split}
        u(x,0) &= 0 \\
        u(0,t) &= u(2,t) = 0 
    \end{split}
    \label{eq: 2DfPDE_ICBC}
\end{equation}
and the source term $f(x,t)$ presented as:
\begin{equation}
    f(x,t) = \frac{2}{\Gamma(3-\alpha)}x(2-x)t^{2-\alpha} + t^2x(2-x) + 2t^2
    \label{eq: 2DfPDE_sourceterm}
\end{equation}
This specific equation can be solved analytically with the following form:

\begin{equation}
    u(x,t) = t^2x(2-x)
    \label{eq: 2DfPDE_analytical}
\end{equation}

The same equation can be also considered as a \textbf{\textit{three-dimensional case}}, and written in the dimensions of $(x,y,t)$ as:
\begin{equation}
\begin{split}
    &\frac{\partial^\alpha u(x,y,t)}{\partial t^\alpha} + u(x,y,t) = \frac{\partial^2 u (x,y,t)}{\partial x^2} + \frac{\partial^2 u (x,y,t)}{\partial y^2} + f(x,y,t)  \\ &0\leq x,y\leq 2, 0\leq t, 0<\alpha<1
    \label{eq: 3DfPDE}
\end{split}
\end{equation}
where with initial and boundary conditions of:
\begin{equation}
    \begin{split}
        u(x,y,0) &= 0\\
        u(0,y,t) &= u(2,y,t) = t^2y(2-y) \\
        u(x,0,t) &= u(x,2,t) = t^2x(2-x)
    \end{split}
    \label{eq: 3DfPDE_ICBC}
\end{equation}
and the source term of:
\begin{equation}
    f(x,y,t) = \frac{2}{\Gamma(3-\alpha)}x(2-x)t^{2-\alpha}[x(2-x) + y(2-y)] + t^2[x(2-x) + y(2-y)] + 4t^2
    \label{eq: 3DfPDE_sourceterm}
\end{equation}
the analytical solution shall be obtained in the form of:
\begin{equation}
    u(x,y,t) = t^2[x(2-x)+y(2-y)]
    \label{eq: 3DfPDE_analytical}
\end{equation}

\subsection{Caputo solution of fractional derivatives}\label{frac_deriv}
The fractional derivative of a function $f(t)$ with respect to time in the Caputo formalism is defined as: 
\begin{equation}
    D^\alpha_t f(t) = \frac{1}{\Gamma(1-\alpha)}\int_0 ^t\frac{f'(x)}{(t-x)^\alpha}dx
    \label{eq: caputo}
\end{equation}
where $\alpha$ and $\Gamma(.)$ are fractional derivative order $(0 \leq \alpha \leq 1)$, and gamma function explained as $\Gamma(x) = \int_0^\infty t^{x-1}e^{-t}dt$, respectively \cite{frac_deriv}. Various numerical methods have been developed to alleviate the challenges stemmed from the computation of the integral term among which, two shall be used in this study. In both methods the interval $[0,t]$ is spaced uniformly with a spacing of $h$, where $t_n = nh, n = 0,1,...,N$.

The first approach  developed by Diethelm et al. \cite{Diethelm_caputo} is based on a finite difference method in which \cref{eq: caputo} is approximated at point $t_r$ through the following equation:
\begin{equation}
    D^\alpha_t f(t_r) = \frac{1}{h^\alpha \Gamma(2-\alpha)}\sum_{n=0}^{n_r}a_{n,n_r}(f_{n_r-n}-f_0)
    \label{eq: Diethelm_caputo}
\end{equation}
where $t_r = n_rh$, and the coefficient $a_{n,n_r}$ is derived as:
    \begin{equation}
    a_{n,n_r} = 
\begin{cases}
      1, & \text{if }n=0 \\
      (n+1)^{1-\alpha}-2n^{1-\alpha}+(n-1)^{1-\alpha}, & \text{if }0<n<n_r \\
      (1-\alpha)n_r^{-\alpha}-n_r^{1-\alpha}+(n_r-1)^{1-\alpha}, & \text{if }n=n_r
    \end{cases}
    \label{eq: anN}
\end{equation}
In the second approach introduced as L1 approximation \cite{L1_caputo}, \cref{eq: caputo} is computed as: 
\begin{equation}
    D^\alpha_t f(t_r) = h^{-\alpha}[b_0f(t_r)-b_{r-1}f(0)+\sum_{n=1}^{r-1}(b_n-b_{n-1})f(t_{r-n})]
    \label{eq: L1_caputo}
\end{equation}
where the coefficient $b_n$ is given as:
\begin{equation}
    b_n = \frac{(n+1)^{1-\alpha}-n^{1-\alpha}}{\Gamma(2-\alpha)}, \quad n = 0,1,...,N-1
    \label{eq: b_n}
\end{equation}
It is worth mentioning that order of accuracy in both methods is $2-\alpha$. Throughout the study, the presented approaches are referred to as Diethelm and L1 method, respectively.
\subsection{Fractional derivatives in neural networks}\label{NN_frac}
\begin{figure}
    \centering
    \includegraphics{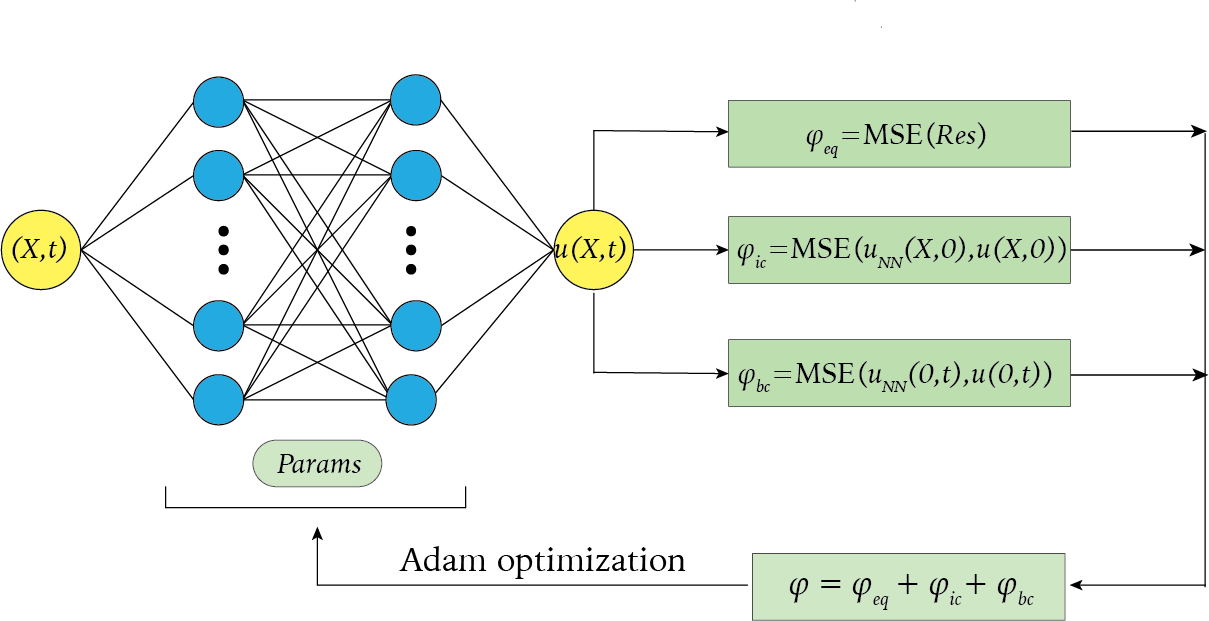}
    \caption{Schematic illustration of the neural network's architecture. The variable $X$ is: nonexistent for fODEs, $x$ for the two-dimensional fPDE, and $(x,y)$ for the three-dimensional fPDE cases.}
    \label{fig:NN}
\end{figure}

A neural network is employed to solve \cref{eq: fODE}, \cref{eq: 2DfPDE} and \cref{eq: 3DfPDE} consisting multiple layers, each containing several neurons. A schematic view of the neural network's architecture is shown in figure \ref{fig:NN}.

\begin{figure}[!htb]
  \includegraphics[width=\linewidth]{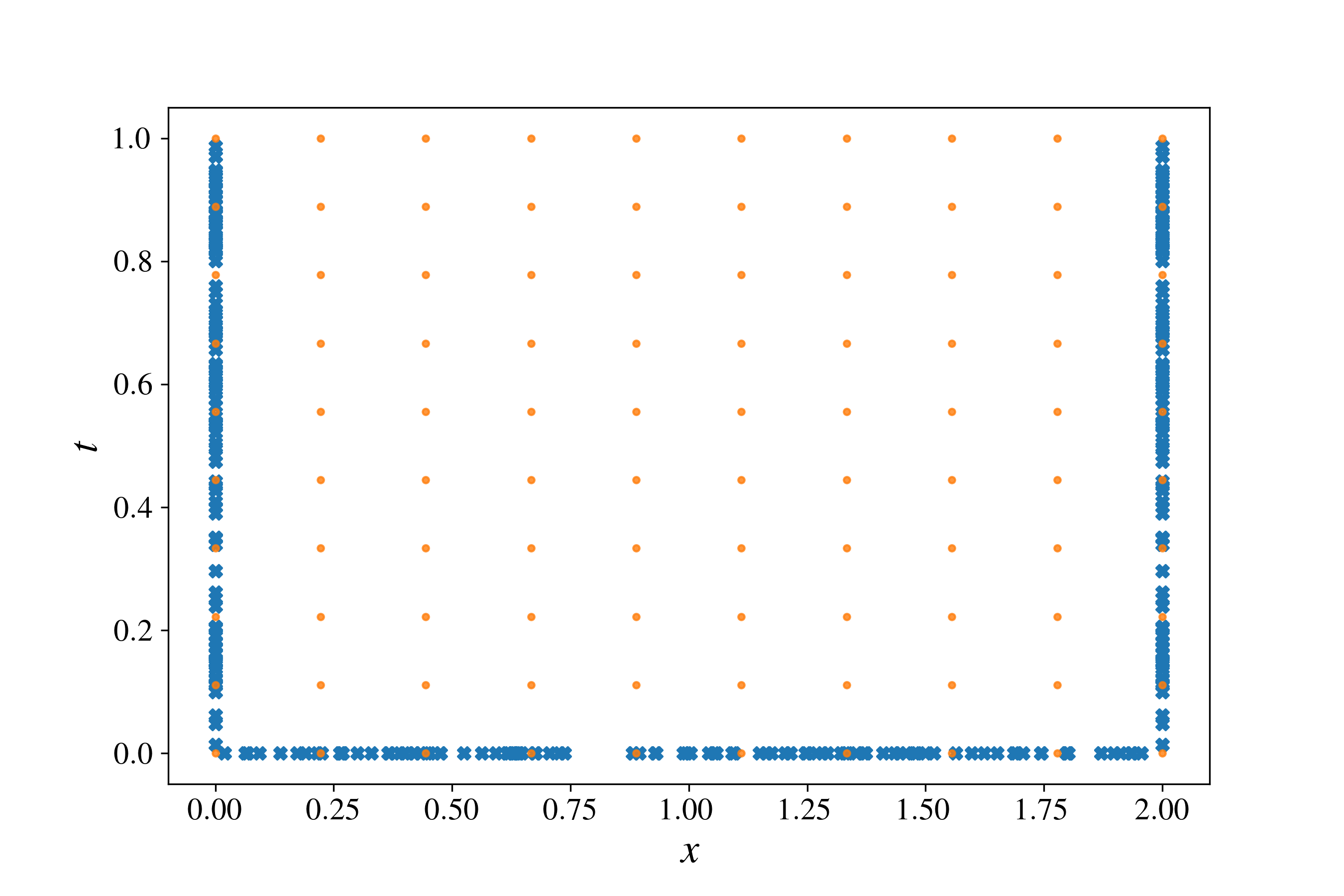}
  \caption{An example of the collocation points defined in the \textit{2D} case with 10 points in $t$, 10 points in $x$, and 100 boundary points and 100 initial points. The orange circles represent the collocation points and the blue crosses represent the boundary data.}\label{fig:collocation_points}
\end{figure}

Firstly, the variable $t$ in the fODE, $(x,t)$ in the 2D, and $(x,y,t)$ in the 3D case are spaced uniformly in their corresponding intervals, giving a group of collocation points in each direction. Details on the number of collocation points are tabulated in \cref{table: hyperparameters}. Upon proceeding, the groups of collocation points are combined, thus giving a grid as presented in figure \ref{fig:collocation_points}. The set of collocation points are given through a tensor as an input to the neural network.

As demonstrated in figure \ref{fig:NN}, in order to train the neural network for solving the equations, neuron weight parameters should be determined by means of minimizing the loss given as:
\begin{equation}
    \phi = \phi_{eq} + \phi_{ic} + \phi_{bc}
    \label{eq: loss_function}
\end{equation}
In the above equation, $\phi_{eq}$, $\phi_{ic}$ and $\phi_{bc}$ are equation loss, loss due to the initial condition, and the loss associated with boundary condition, respectively, each of which corresponding to a specific category of collocation points.

Equation loss is correlated to the equation residual. The residual of \cref{eq: fODE}, \cref{eq: 2DfPDE} and \cref{eq: 3DfPDE} are given as:

\begin{equation}
    \begin{split}
        Res(t) &= \frac{\partial^\alpha u_{NN}(t)}{\partial t^\alpha} + u_{NN}(t) - t^2 - \frac{8}{3\sqrt{\pi}}t^{\frac{3}{2}}\\
        Res(x,t) &= \frac{\partial^\alpha u_{NN}(x,t)}{\partial t^\alpha} + u_{NN}(x,t) - \frac{\partial^2 u_{NN}(x,t)}{\partial x^2} - f(x,t) \\
        Res(x,y,t) &= \frac{\partial^\alpha u_{NN}(x,y,t)}{\partial t^\alpha} + u_{NN}(x,y,t) - (\frac{\partial^2 u_{NN}(x,y,t)}{\partial x^2} + \frac{\partial^2 u_{NN}(x,y,t)}{\partial y^2})\\ &- f(x,y,t)
    \end{split}
    \label{eq:residuals}
\end{equation}
where $u_{NN}$ represents the neural network output. As indicated in \cref{eq:residuals}, fractional derivative methods are utilized here. Consequently, their functionalities can be evaluated through $\phi_{eq}$, and actions shall be measured to improve the training, in case of encountering high values of $\phi_{eq}$. Residual is computed in all sets of the collocation points and is then employed in computing $\phi_{eq}$ through the following mean-squared error (MSE)
\begin{equation}
    \phi_{eq} = MSE(Res) = \frac{1}{N_{eq}}\sum_{n=1}^{N_{eq}}Res^2(X_n,t_n)
    \label{eq:MSE(residual)}
\end{equation}
where $N_{eq}$ is the number of collocation points for residual and $X_n$ is nonexistent, $x_n$ and $(x_n, y_n)$ in fODE, two-dimensional and three-dimensional fPDE cases, respectively.

The losses due to initial and boundary conditions follow similar procedures as they present the discrepancy between the neural network output in initial time/boundaries and the initial/boundary conditions.
\begin{equation}
    \begin{split}
        \phi_{ic} = MSE(u_{NN}(X,0)-u(X,0)) = \frac{1}{N_{ic}}\sum_{n=1}^{N_{ic}}(u_{NN}(X_n,0)-u(X_n,0))^2\\
        \phi_{bc} = MSE(u_{NN}(0,t)-u(0,t))= \frac{1}{N_{bc}}\sum_{n=1}^{N_{bc}}(u_{NN}(0,t_n)-u(0,t_n))^2
    \end{split}
    \label{eq:MSE_ICBC}
\end{equation}
where $N_{ic}$ and $N_{bc}$ are the number of collocation points set for computing initial condition and boundary condition loss. Details on the neural network hyperparameters, including number of layers, neurons per layer, etc., are tabulated in \cref{table: hyperparameters} for each of the cases studied here. In all the cases, learning rate is chosen in a step-wise form to decrease as the training proceeds. This approach leverages high training speed in the preliminary iterations while in the subsequent iterations, neural network avoids missing the optimum point since it has a small learning rate.

\begin{table}[h]
\caption{Hyperparameters used in training neural network for solving fODE (\cref{eq: fODE}), two-dimensional fPDE (\cref{eq: 2DfPDE}) and three-dimensional fPDE (\cref{eq: 3DfPDE})}\label{table: hyperparameters}%
\begin{small}
\begin{tabularx}{\textwidth}{XXXX}
\toprule
Hyperparameter & Equation \ref{eq: fODE} & Equation \ref{eq: 2DfPDE} & Equation \ref{eq: 3DfPDE}\\
\midrule
Fractional derivative order ($\alpha$) & 0.5 & 0.5 & 0.5 \\
Iterations of learning rate change & 200,1000 & 2000,5000 & 2000,5000 \\ 
Learning rate   & 0.01,0.001,0.0005 & 0.01,0.005,0.001  & 0.01,0.005,0.001\\
Number of layers    & 3 & 4 & 4\\
Number of neurons    & 10 & 20 & 20\\
Equation domain points\\ ($t$)/$(x,t)$/$(x,y,t)$    & 30 & (10,10) & (5,5,5)\\
Initial condition points  & 30 & 100 & 5\\
Boundary condition points   & - & 100 & 25\\
\bottomrule
\end{tabularx}
\end{small}
\end{table}

\section{Results and Discussion}\label{sec: result}

The ultimate goal of the present work is to: (1) investigate the compromise between accuracy and computational costs, and seek plausible approaches to improve the neural network's results with imposing minimum burden on computational cost, and (2) provide a detailed guideline on devising fPINNs with high levels of efficiency and accuracy. To this end, the results for each of the defined equations in \cref{sec: method} are presented and discussed in the same order. 

\begin{figure}
    \centering
    \includegraphics[width=0.7\textwidth]{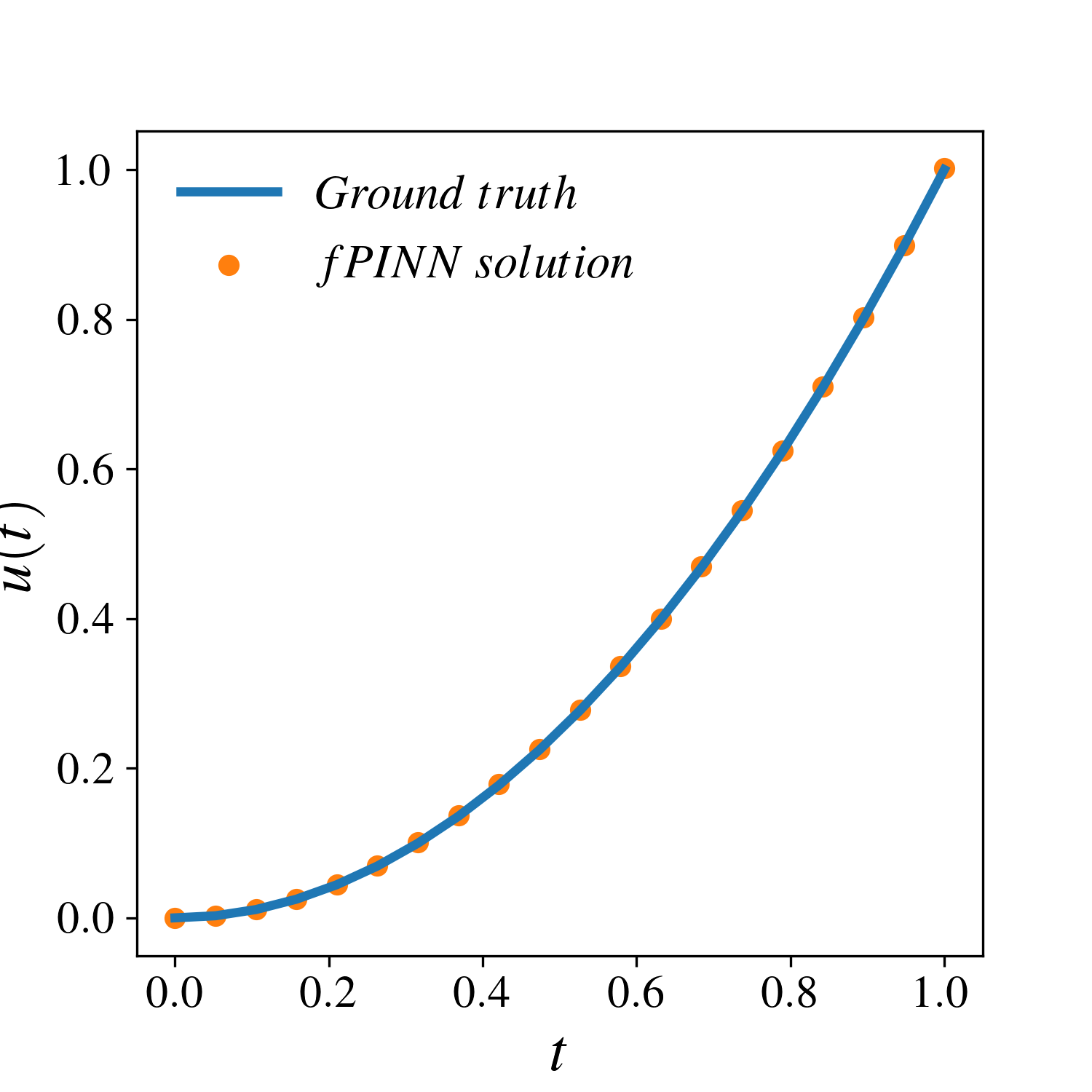}
    \caption{A comparison between fPINN results using the Diethelm method and the analytical solution of \cref{eq: fODE}.}
    \label{fig:fODE}
\end{figure}

Employing the hyperparamters tabulated in \cref{table: hyperparameters}, fPINN solves \cref{eq: fODE}, using Diethelm method within a few minutes. As presented in figure \ref{fig:fODE}, neural network accurately predicts the result throughout the whole domain. Expectedly, for a simple fODE, data-driven solutions and the ground truth are equivalent.

With increasing the level of complexity to a two-dimensional fPDE, fPINN is next applied to the \cref{eq: 2DfPDE}. The details of the hyperparameters used can be found in \cref{table: hyperparameters}. In training the neural network and solving \cref{eq: 2DfPDE}, both Diethelm and L1 approaches are employed to give an insight into their accuracy for this specific problem. The neural network is trained for 10 minutes. 

\begin{figure}
    \centering
    \includegraphics[width=1.\textwidth]{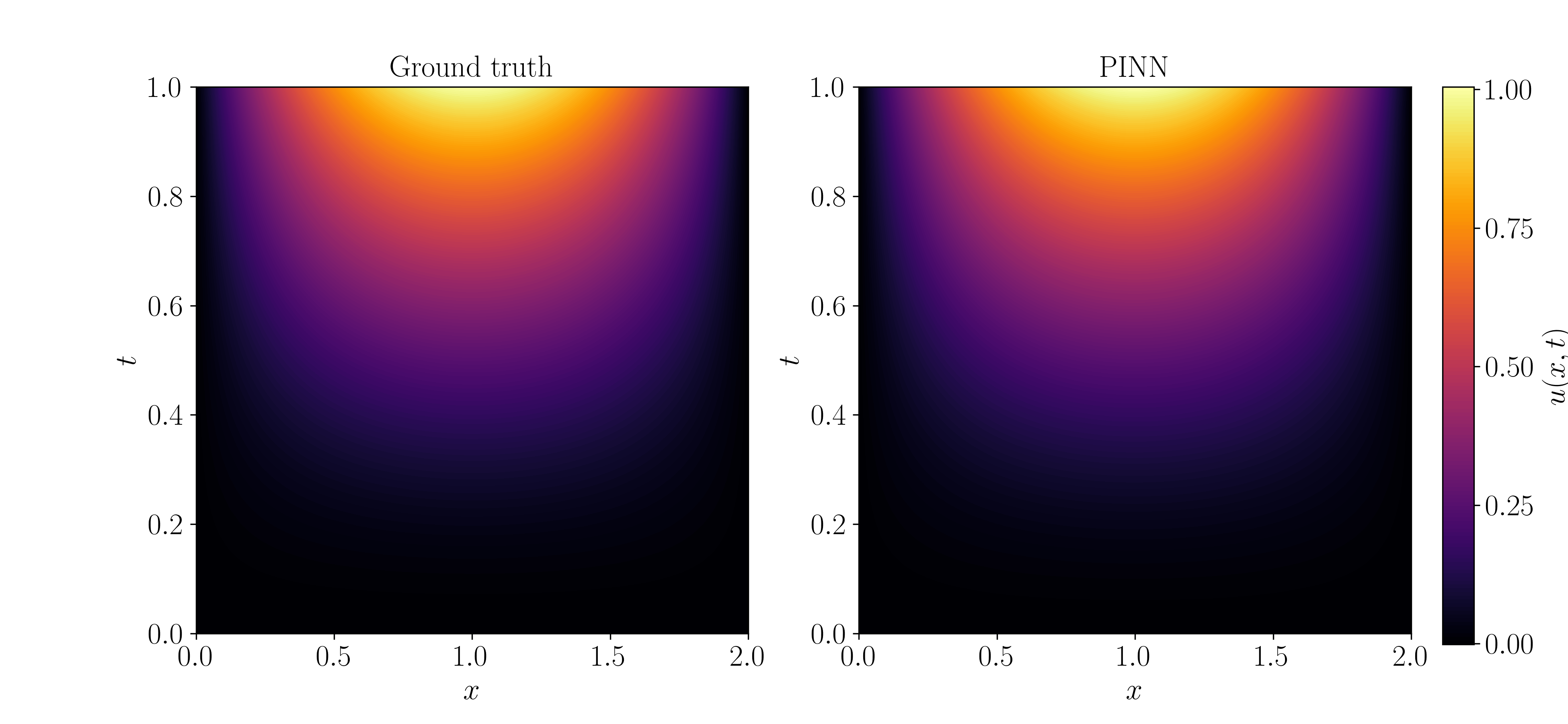}
    \caption{A comparison between fPINN result and the analytical solution of \cref{eq: 2DfPDE} in a given domain. The corresponding hyperparameters are presented in \cref{table: hyperparameters} and the exact value of the function $u$ is given in the colorbar.}
    \label{fig:2DfPDE_t10_contour}
\end{figure}

\begin{figure}
    \centering
    \includegraphics[width=1.\textwidth]{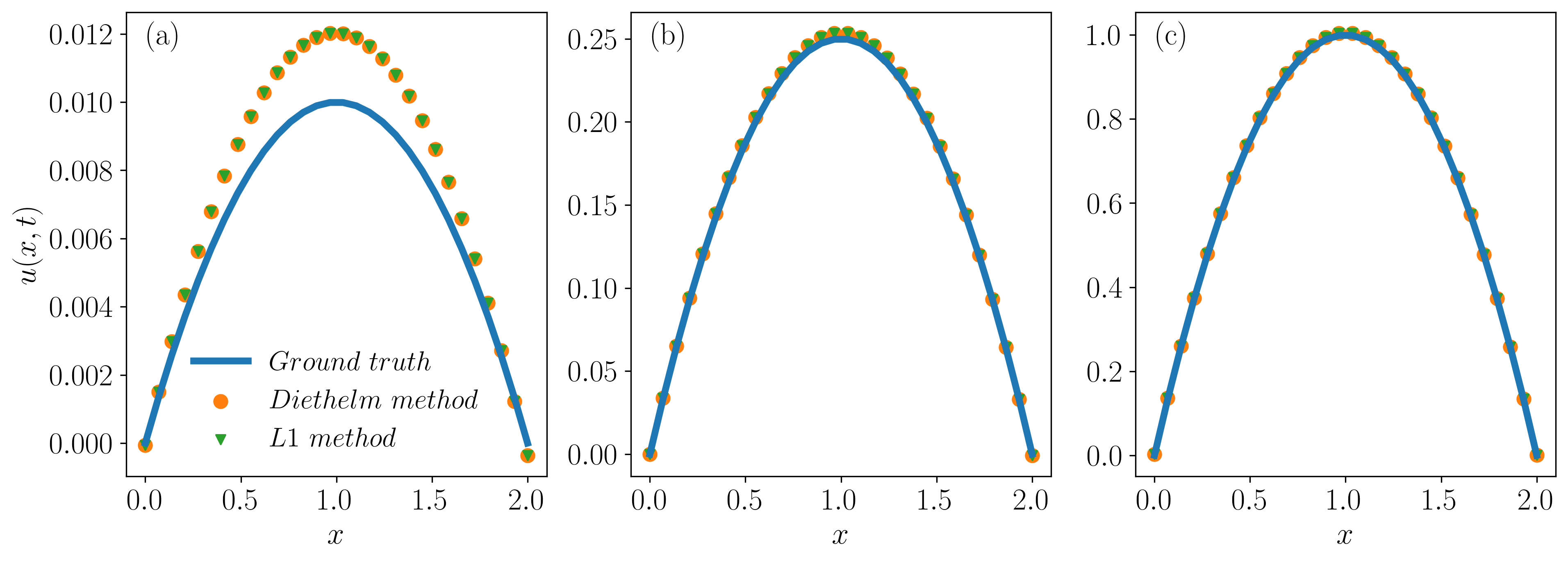}
    \caption{A comparison between fPINN results and analytical solution of \cref{eq: 2DfPDE} at different times of (a) $0.1 s$, (b) $0.5 s$ and (c) $1 s$, using both Diethelm and L1 method.}
    \label{fig:2DfPDE_t10}
\end{figure}

As presented in figure \ref{fig:2DfPDE_t10_contour}, the fPINN results closely follow the ground truth solution across the domain. Nonetheless and to get a closer insight, fPINN solutions are plotted in figure \ref{fig:2DfPDE_t10} in the $x$ direction and at three different time stamps of 0.1\textit{s}, 0.5\textit{s}, and 1\textit{s}. Results in figure \ref{fig:2DfPDE_t10} indicate that while at longer times the data-driven and the numerical solutions are accurately matching, there are slight deviations at the early times for both methods used (Diethelm and L1). This is largely because in both proposed methods (\cref{eq: Diethelm_caputo} and \cref{eq: L1_caputo}), the information of previous points is required for computing fractional derivative at a point. This inherently results in deviations to be observed at initial points, and where the least amount of data (or no data at all) has been observed. As the time progresses and the history of a point is observed with more data, fPINN solutions and the numerical solutions begin to converge.

In the following, we investigate different strategies to improve upon the accuracy of these predictions. Namely, we will interrogate the role of increasing collocation points in $t$ direction, and reducing the time interval with a constant number of collocation points. In the following, both strategies and their effects on results-in terms of accuracy and computational cost-are studied. Since in figure \ref{fig:2DfPDE_t10}, both L1 and Diethelm methods offer identical results, only the results of Diethelm method are presented in the following discussions.

\begin{figure}
    \centering
    \includegraphics[width=1.\textwidth]{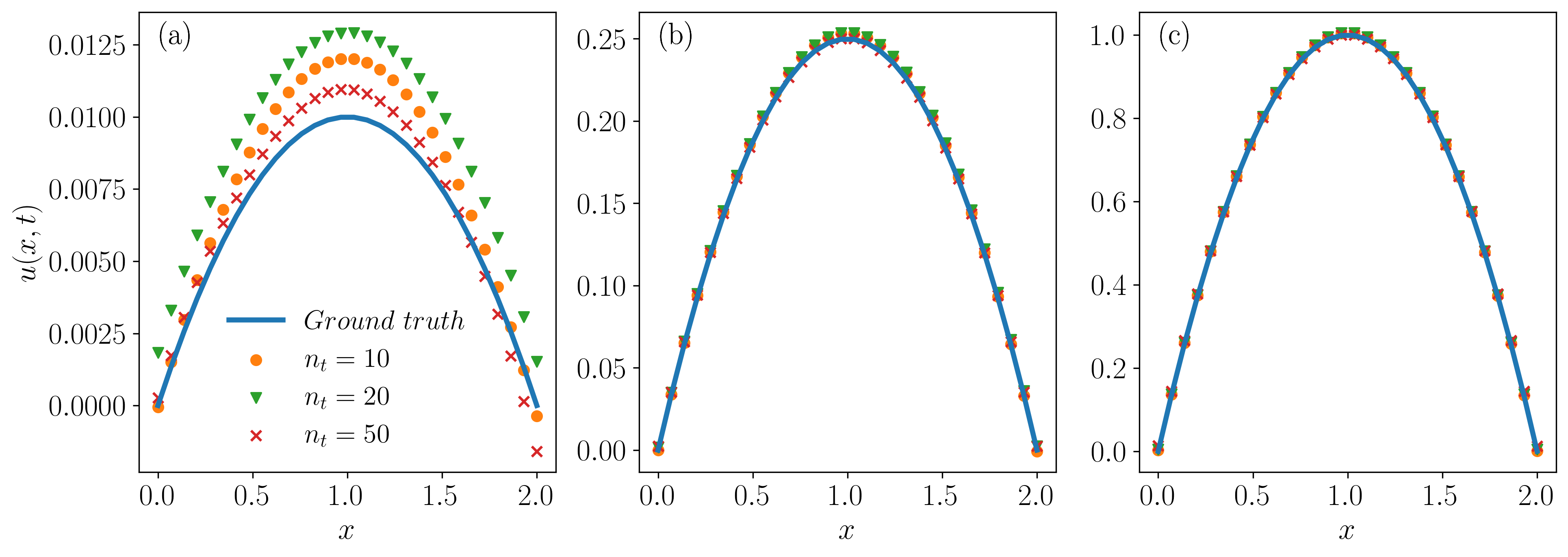}
    \caption{fPINN solution of \cref{eq: 2DfPDE} using different number of collocation point sets at different times of (a) 0.1 \textit{s}, (b) 0.5 \textit{s} and (c) 1 \textit{s}, benchmarked against the ground truth solution. The legend indicates the number of collocation points in the $t$ direction. In the simulations, the number of collocation points is increased equally in both $x$ and $t$ directions. The results are obtained using Diethelm method.}
    \label{fig:2DfPDE_colloc_effect}
\end{figure}

In order to study the role of collocation points on the final results of the fPINN, collocation points are increased in both $x$ and $t$ directions from 10 to 20, and 50, respectively, with run times reported in \cref{table: 2DfPDE_runtimes}. Generally, increasing the number of collocation points plays a vital role in the accuracy of final predictions, as shown in figure \ref{fig:2DfPDE_colloc_effect}. Doubling the number of collocation points in each direction from 10 to 20 along with increasing the run time enhances the accuracy at initial times. Additionally, with a constant run time, adding data points is also effective in increasing the accuracy despite the reduced number of iterations that the neural network goes through.

\begin{table}[h]
\caption{Training specifics for each set of collocation points in the two-dimensional fPDE, using Diethelm method. In all cases, number of collocation points in the $x$ and $t$ directions increase to the same number of points.}\label{table: 2DfPDE_runtimes}%
\begin{small}
\begin{tabularx}{\textwidth}{XXXX}
\toprule
Number of collocation points in time & Run time (minutes) & Number of iterations & Total loss order\\
\midrule
10    & 10 & $9.3\times10^{4}$ & $10^{-6}$\\
20    & 30 & $5.74\times10^{4}$ & $10^{-5}$\\
50    & 30 & $4.4\times10^{3}$ & $10^{-5}$\\
\bottomrule
\end{tabularx}
\end{small}
\end{table}

As shown in figure \ref{fig:2DfPDE_colloc_effect}, increasing the number of collocation points improves upon the accuracy of the fPINN predictions at the early times, and also increases the overall runtime significantly. As such, and in scenarios where one is not as much interested in the results at the initial time stamps, running fPINN with smaller number of collocation points may be of interest. However, if the accuracy of solution in the initial time is of particular interest, reducing time window can be advantageous. With this arrangement collocation points are spaced in a smaller interval, giving smaller time steps and consequently lower error in computing fractional derivative order as the order of accuracy is $2-\alpha$.

\begin{figure}
    \centering
    \includegraphics[width=1.\textwidth]{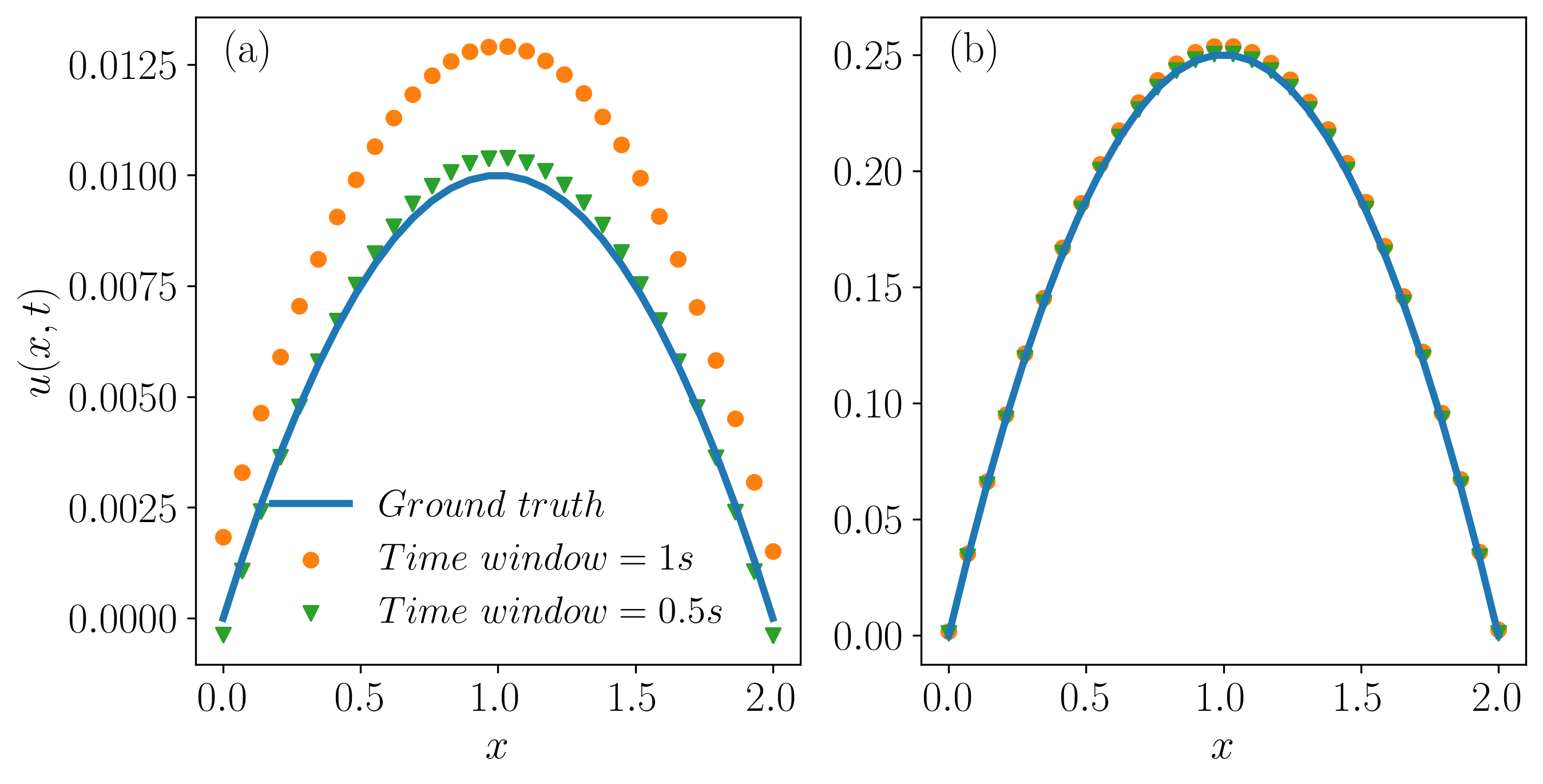}
    \caption{fPINN results of \cref{eq: 2DfPDE} with 20 collocation points in time over a window of 1\textit{s} and 0.5\textit{s} at different timestamps of (a) 0.1 \textit{s}, (b) 0.5 \textit{s}. Neural network results are obtained using Diethelm method.
    }
    \label{fig:2DfPDE_timewindow_effect}
\end{figure}

In figure \ref{fig:2DfPDE_timewindow_effect}, closer agreement between the fPINN prediction and the ground truth solution in panel (a) indicates the effect of reduced time window on improving accuracy in the number of 20 collocation points in $t$ and 30 minutes of run time.

\begin{figure}
    \centering
    \includegraphics[width=1.\textwidth]{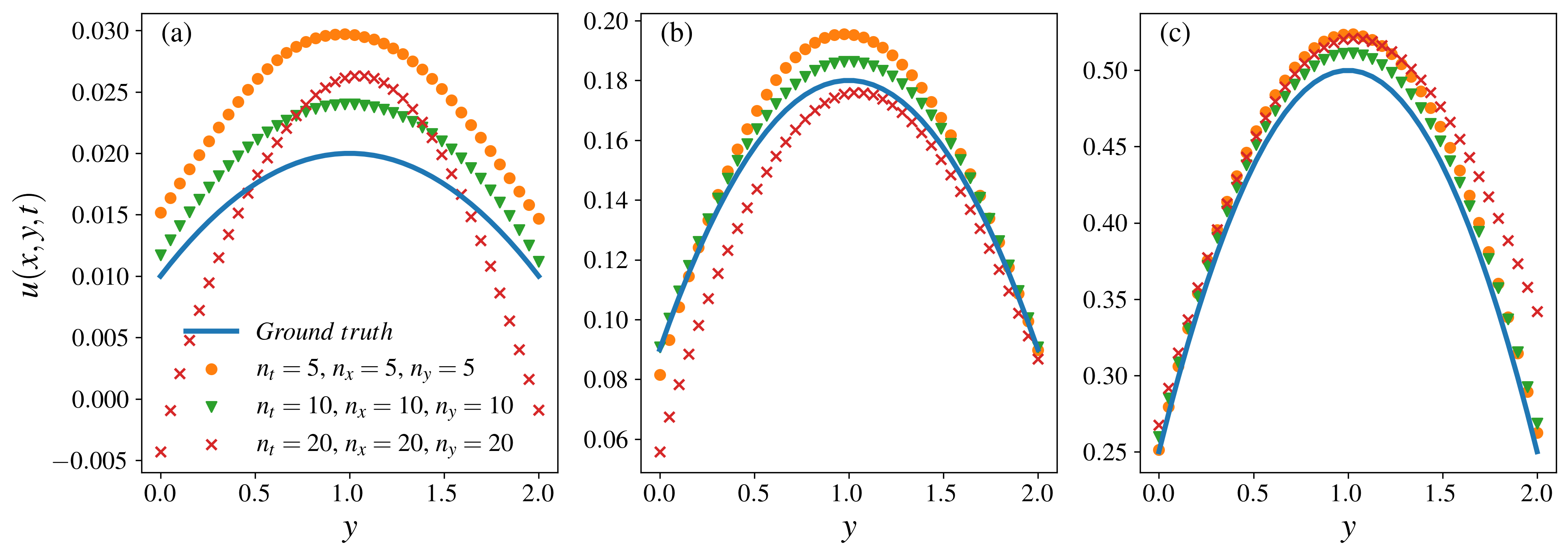}
    \caption{A comparison between fPINN results and the ground truth solution of \cref{eq: 3DfPDE} at different times of (a) 0.1 \textit{s}, (b) 0.3 \textit{s}, and (c) 0.5 \textit{s}, using different sets of collocation points in $t$, $x$ and $y$ directions, as indicated in the legend. The results are plotted at the half-length of domain ($x = 1$).}
    \label{fig:3DfPDE_EvenSets}
\end{figure}

Next, fPINNs are trained to solve a three dimensional fPDE, given in \cref{eq: 3DfPDE}, utilizing hyperparameters outlined in \cref{table: hyperparameters} with the exception of time window where time varies from zero to 0.5 \textit{s}, with a run time of 30 minutes to have an estimation of the necessary number of collocation points. For three sets of (5,5,5), (10,10,10) and (20,20,20) collocation points in 
the $x$, $y$, and $t$ directions respectively, the fPINN results in the the domain's half-length ($x = 1$) are presented in figure \ref{fig:3DfPDE_EvenSets}. 

\begin{figure}
    \centering
    \includegraphics[width=1.\textwidth]{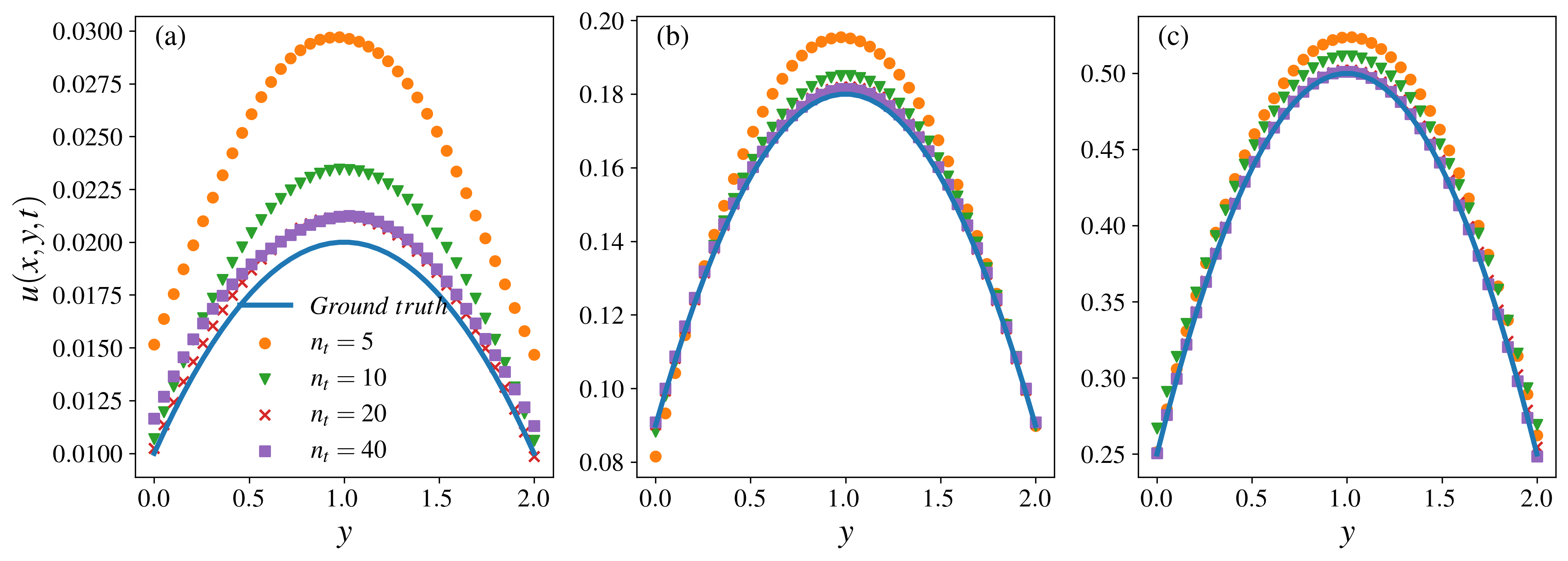}
    \caption{fPINN solutions of \cref{eq: 3DfPDE}
    benchmarked against the ground truth solution for different volumes of collocation points used during training at different times of (a) 0.1 \textit{s}, (b) 0.3 \textit{s} and (c) 0.5 \textit{s}. The number of points in $(x,y)$ are kept constant at (5,5) and the number of points in the $t$ direction are varied from 5 to 40, given in the legend.}
    \label{fig:3DfPDE_time_colloc_effect}
\end{figure}

The results in figure \ref{fig:3DfPDE_EvenSets} clearly show that increasing the number of data points for the fPINN solution can have a non-monotonic effect. Initially, increasing the number of collocation points in each parameter from 5 to 10 enhances the accuracy, specially for the initial times (the long time solutions remain very close for these two cases). Nonetheless, increasing the number of collocation points further to 20 points, results in an incomplete training over the time allocated, and thus yielding an erroneous solution which at early times is largely deviating from all other solutions, and at longer times does not seem to improve significantly. It is worth noting that in all three cases the neural network's results are not as accurate as the fODE or the two dimensional fPDE solutions. In contrast to the previous scenario, increasing the number of collocation points does not generally result in better accuracy, unless the overall run time is also significantly increased. Naturally and with no bounds to the run time, more data points will result in better predictions; however, in the range of studying parameters and equations here, increasing the collocation points for the three dimensional case inevitably makes the fPINN impractical. On the other hand, and having three independent variables of $(x,y,t)$, raises the question of whether changing the distribution of collocation points in one dimension preferentially, will result in a better performance. In other words, what if keeping an overall limited number of data points (fewer than 8000 corresponding to the (20,20,20) case), better results can be achieved merely by a different distribution of those points in each direction? Setting the case of (5,5,5) collocation set as the smallest number of data points, the number of collocation points are increased in the $t$ direction, keeping an overall run time of 30 minutes constant and the results are presented in figure \ref{fig:3DfPDE_time_colloc_effect}. Evidently, increasing the collocation points in one direction (time) from 5 to 20, as opposed to increasing it in all directions improves the accuracy of fPINN solution to an acceptable extent in which the fPINN solution only slightly deviates from the ground truth solution at the initial time. At longer times, these conditions consistently provide accurate solutions. Additionally, increasing the collocation points from 20 to 40 does not seem to change the overall accuracy of solution significantly.

\begin{figure}
    \centering
    \includegraphics[width=1.\textwidth]{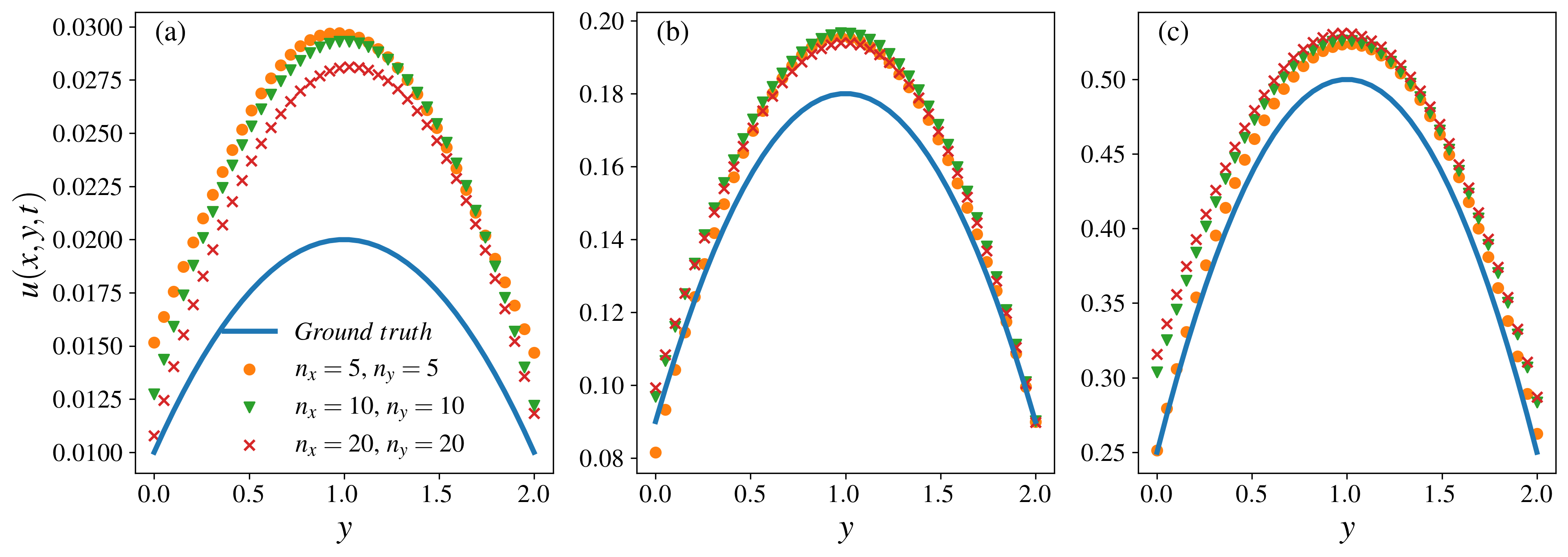}
    \caption{fPINN solutions of \cref{eq: 3DfPDE} using different numbers of collocation points in $x$ and $y$ direction benchmarked against the ground truth solution 
     at timestamps of (a) 0.1 \textit{s}, (b) 0.3 \textit{s} and (c) 0.5 \textit{s}. The number of points in $t$ direction is 5.}
    \label{fig:3DfPDE_xy_colloc_effect}
\end{figure}

The same effect was also studied for the other two variables and by keeping the number of points in the $t$ direction constant while changing them in the $(x,y)$ directions. Interestingly, results in figure \ref{fig:3DfPDE_xy_colloc_effect} clearly show that increasing the number of data points in the $(x,y)$ directions does not significantly improve the results. As such, these results suggest that in cases where collection of data is possible for training purposes, priority should be given to adding sensory data to time rather than space.

\begin{figure}
    \centering
    \includegraphics[width=1.\textwidth]{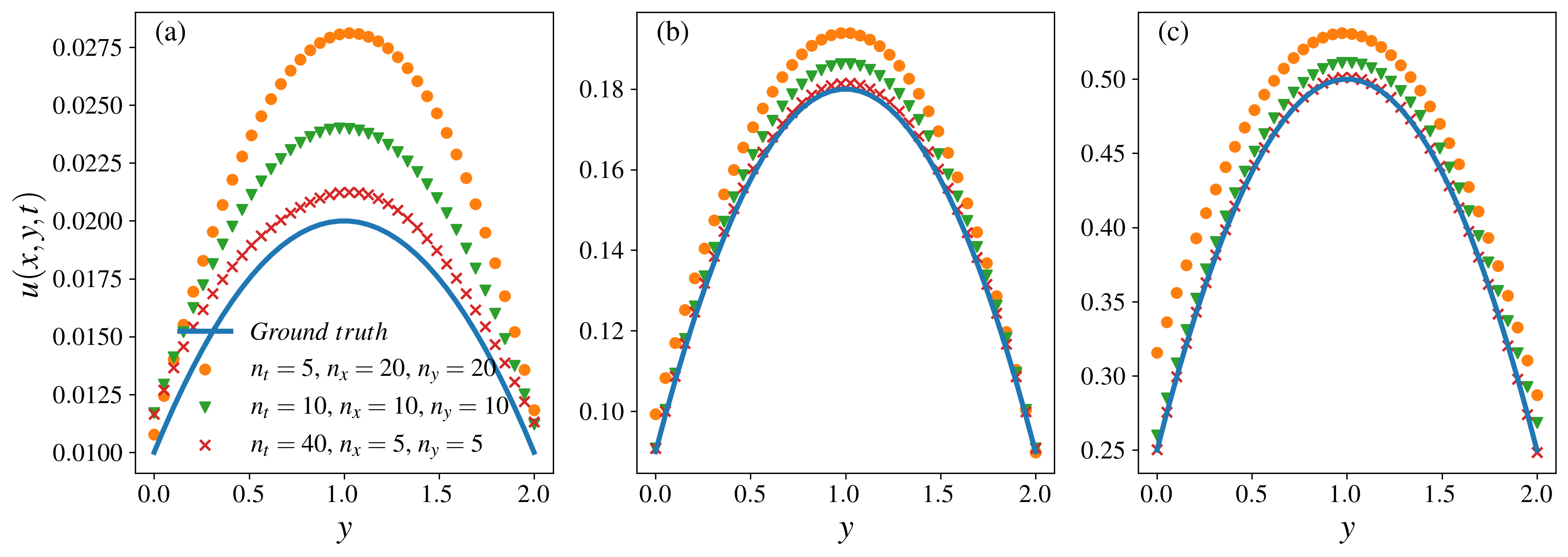}
    \caption{
    A comparison between fPINN solutions of \cref{eq: 3DfPDE} and ground truth solution for three different distributions of 1000 collocation points in times of (a) 0.1 \textit{s}, (b) 0.3 \textit{s} and (c) 0.5 \textit{s}.}
    \label{fig:3DfPDE_different_distributions}
\end{figure}

In order to give a side-by-side comparison of different data sets, the best results in each of the previous figures are brought into comparison in figure \ref{fig:3DfPDE_different_distributions}, where a total number of 1000 data points are distributed in three distinct ways: Evenly, preferentially adding points points in $t$, and preferentially adding points points in $(x,y)$. Evidently, having additional data points in $t$ is the most effective route to increasing the overall solution accuracy. Tensorflow uses GradientTape calculating the integer derivatives which has a high accuracy in terms of discretization. On the contrary, as previously mentioned, the two proposed methods (Diethelm and L1) for deriving fractional derivatives have an order of accuracy $2-\alpha$ which is close to 2, in the lowest value of $\alpha$. Hence, the result is more impacted by the the data in the direction with lower accuracy, which is $t$ in this problem.

\begin{figure}
    \centering
    \includegraphics[width=1.\textwidth]{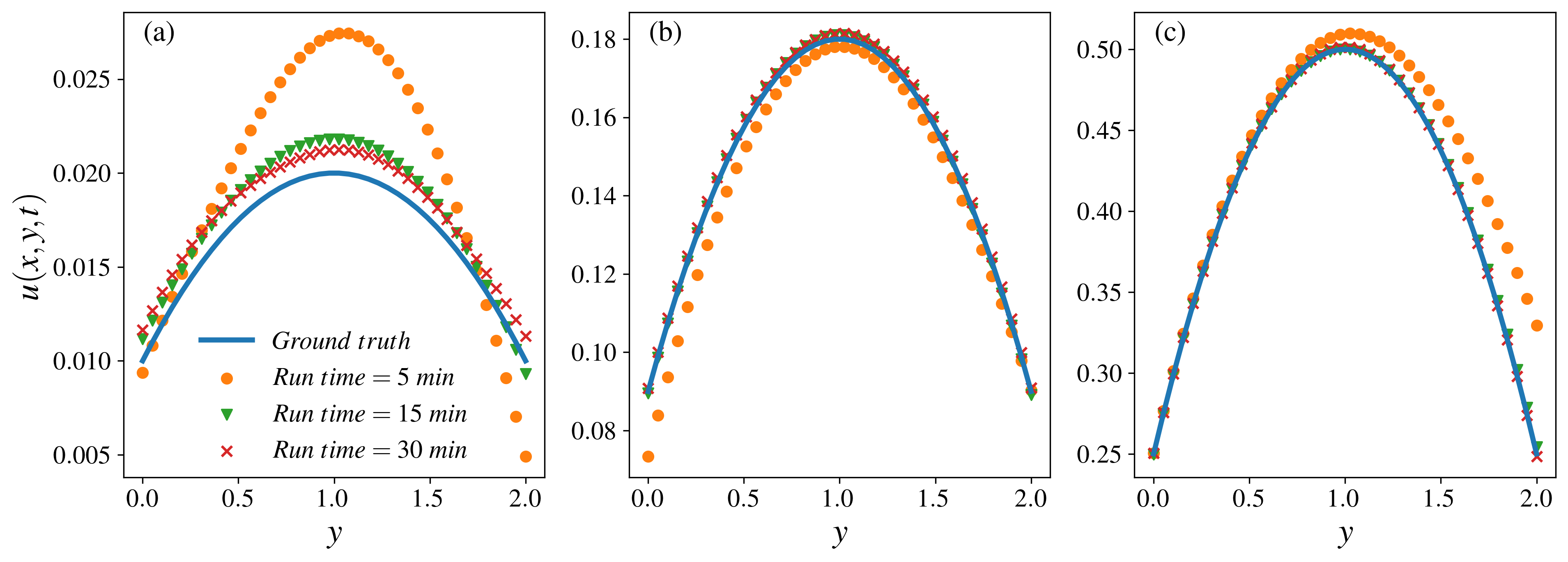}
    \caption{fPINN solutions of \cref{eq: 3DfPDE} benchmarked against the ground truth solution for the set of (5,5,40) collocation points in $x$, $y$ and $t$ directions respectively, with 5, 15 and 30 minutes of run time in times of (a) 0.1 \textit{s}, (b) 0.3 \textit{s} and (c) 0.5 \textit{s}.}
    \label{fig:3DfPDE_different_runtimes}
\end{figure}

Note that while the fPINN solution with training performed using a (5,5,40) data set offers the most accurate solution, all of these solutions are presented using an overall run time of 30 minutes. One should also consider the computational cost of training the fPINN. Our results in figure \ref{fig:3DfPDE_different_runtimes} suggest that for the same set of collocation points, fPINN predictions stay virtually unchanged with half the run time (15 minutes) but begin to deviate largely from the ground truth solutions at shorter run times of 5 minutes. It is worth mentioning that if solutions at longer times are of particular interest, even small training times of 5 minutes may result in sufficiently accurate solutions.

\begin{figure}
    \centering
    \includegraphics[width=1.\textwidth]{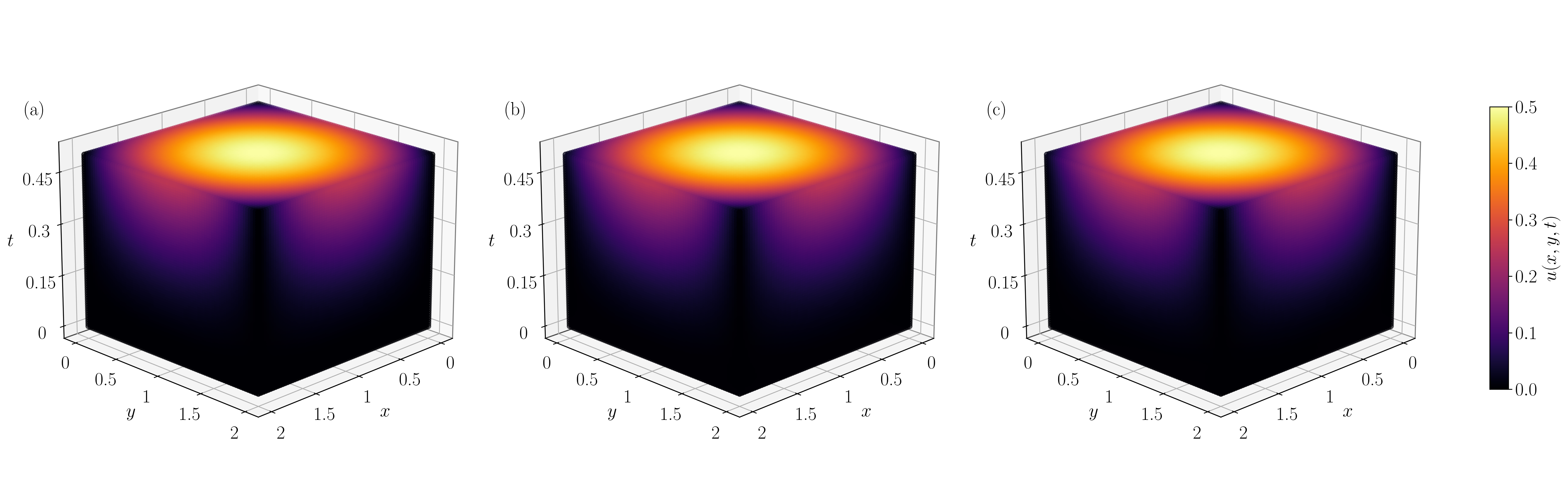}
    \caption{A comparison between (a) ground truth solution and the fPINN solution of \cref{eq: 3DfPDE} using (b) Diethelm and (c) L1 method
    over the entire domain, using a data set of size (5,5,40) in $x$, $y$, and $t$ direction, respectively with an overall run time of 30 minutes.}
    \label{fig:3DfPDE_3Dplot}
\end{figure}

To provide an overall comparison for the solutions with the optimal parameters, fPINN solution over the entire domain is benchmarked against the ground truth solution in figure \ref{fig:3DfPDE_3Dplot}.

In all of the presented results, Diethelm was the method of choice for computing the fractional derivative. Using the same number of collocation points and overall run time, fPINN was also trained using the L1 method and its solution for \cref{eq: 3DfPDE} over the entire domain was benchmarked against the ground truth solution and compared with the result achieved by Diethelm method in \cref{table: 3DfPDE_L1_vs_Caputo} and figure \ref{fig:3DfPDE_3Dplot}. While both methods are efficient at longer times and very accurate, at the initial times (where most errors are observed consistently), L1 offers a slightly more accurate solution.

\begin{table}[h]
\caption{Training specifics of fPINN with a collocation set of (5,5,40) in $x$, $y$ and $t$ direction and a total run time of 30 minutes using L1 and Diethelm method.}\label{table: 3DfPDE_L1_vs_Caputo}%
\begin{small}
\begin{tabularx}{\textwidth}{XXX}
\toprule
Derivative method & Number of iterations & Total loss\\
\midrule
Diethelm    & 25200 & $1.53\times10^{-5}$\\
L1   & 27100 & $6.57\times10^{-6}$\\
\bottomrule
\end{tabularx}
\end{small}
\end{table}

\begin{figure}
    \centering
    \includegraphics[width=1.\textwidth]{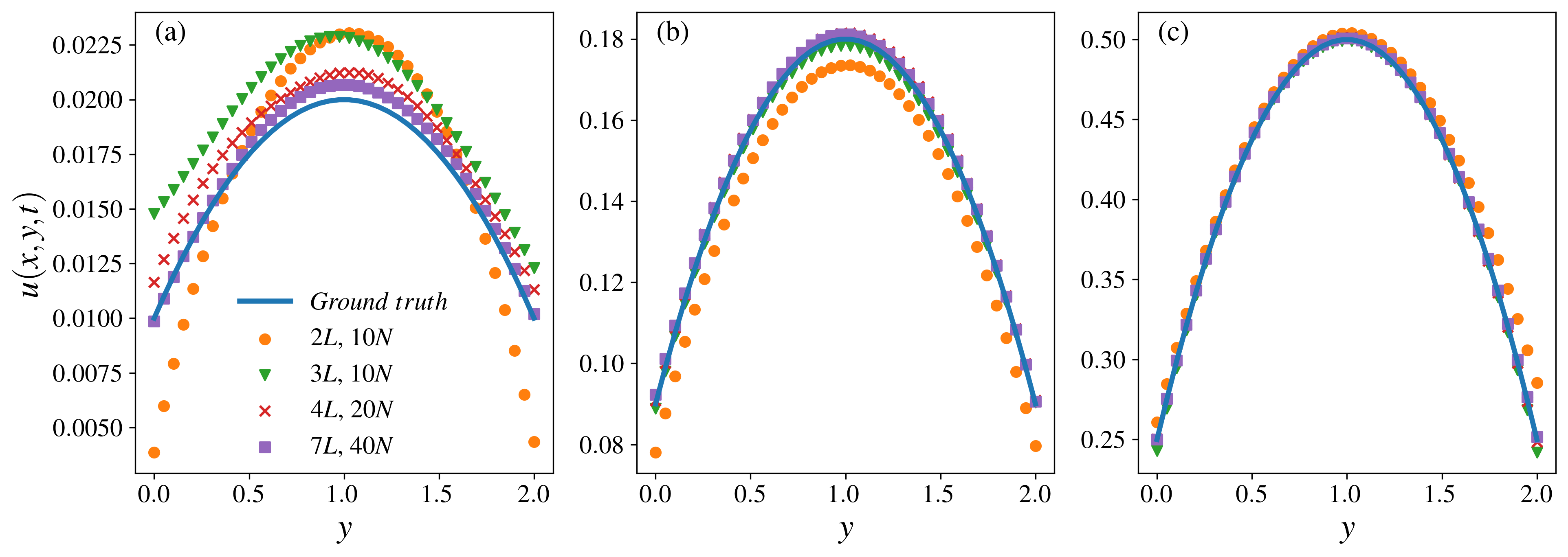}
    \caption{A comparison between results of the fPINN with different hyperparameters (details on \cref{table: hyperparameters}). Solutions of \cref{eq: 3DfPDE}
    are presented at different times of (a) 0.1 \textit{s}, (b) 0.3 \textit{s}, and (c) 0.5 \textit{s} for three different networks of 2, 3, 4 and 7 layers with 10, 10, 20 and 40 neurons per layer, respectively.
}
    \label{fig:3DfPDE_NNarch}
\end{figure}

Lastly, we study the role of fPINN's hyperparameters on its solution accuracy for the most complex cases presented. A wider and deeper neural network containing numerous layers and neurons has a greater number of learning parameters, allowing it to tackle complex problems. However, this comes at the cost of increased computational demand due to the excessive number of learning parameters involved. Carrying the (5,5,40) collocation set and 30 minutes run time from before as the benchmarking case, the training is repeated with larger and smaller neural network architectures to investigate whether the same results can be achieved with a smaller network to reduce the run time. As evident from figure \ref{fig:3DfPDE_NNarch}, shallower networks result in a reduced accuracy, while over-extending the network over four layers does not necessarily result in a significant improvement on the performance. All networks used at long times eventually converge and provide an excellent solution, with the deeper networks yielding marginally better solutions at earlier time stamps. It is of course important to mention that finding the optimal architecture is an iterative process and is highly case dependent.

\section{Conclusion}\label{sec:conclusion}
This study delved into exploring the capability of fPINNs in solving fractional differential equations with different levels of complexity. Our results showed that while fractional ODEs such as one presented in \cref{eq: fODE}, can be precisely solved using fPINNs, solving fractional PDEs (similar to ones presented in \cref{eq: 2DfPDE} and \cref{eq: 3DfPDE}) pose different challenges with respect to the training process and are also associated with varying levels of inaccuracies. In particular, discrepancies are observed between the fPINN solutions and the ground truth solutions at the initial time steps due to the history-based numerical algorithms that were introduced in \cref{sec: method}. To alleviate these issues, one can increase the number of collocation points used during the fPINN training. Nonetheless, employing additional collocation points entails heavier computational cost, creating a trade off between accuracy and runtime. We also found that adding collocation points in the time domain is always more advantageous compared to adding points in the space domain. As such, and when possible, adding sensors to collect data in the time domain should be prioritized over data collection in the space domain. If one is particularly interested in the fPINN solutions at initial times, and considering lack of enough history in those times to yield accurate solutions, reducing the time window with constant collocation points can improve the overall performance of fPINNs. Of course, feasible approaches to enhancing the fPINN performance or reducing the computational cost are not limited to the size of collocation points or the time window. Other parameters such as the size of the fPINN architecture can also play an important role. Once such practical details have been studied and carefully optimized for a given problem, fPINNs present extremely powerful data-driven solutions to a wide range of fractional differential equations.

\section*{Acknowledgements}
The authors acknowledge the support from the National Science Foundation's DMREF program through Award \#2118962.
\label{acknowledgements}

\section*{Declarations}

The authors declare that there is no conflict of interest.





\bibliographystyle{elsarticle-num}
\bibliography{bibliography}
\end{document}